\documentclass[11pt]{article}
\usepackage{amsmath}
\usepackage{amsfonts}
\usepackage[bookmarksnumbered=true]{hyperref}

\newtheorem{theorem}{Theorem}[section]
\newtheorem{lemma}[theorem]{Lemma}
\newtheorem{definition}[theorem]{Definition}
\newtheorem{remark}[theorem]{Remark}

\begin{document}

\title{A Class of Backward Doubly
Stochastic Differential Equations with Discontinuous Coefficients
\thanks{This work is supported by National Natural Science Foundation of China Grant
10771122, Natural Science Foundation of Shandong Province of China
Grant Y2006A08 and National Basic Research Program of China (973
Program, No.2007CB814900)}}
\author{Qingfeng Zhu$^{\rm a}$ and Yufeng Shi$^{\rm b}$
\thanks{Corresponding author, E-mail: yfshi@sdu.edu.cn}\\
{\small $^{\rm a}$ School of Statistics and Mathematics, Shandong
University of Finance},\\
{\small Jinan 250014, China}\\
{\small$^{\rm b}$School of Mathematics, Shandong University, Jinan 250100, China}}
\maketitle

\begin{abstract}In this work the existence of solutions of one-dimensional backward
doubly stochastic differential equations (BDSDEs in short) where the
coefficient is left-Lipschitz in $y$ (may be discontinuous) and
Lipschitz in $z$ is studied. Also, the associated comparison theorem
is obtained.\\
\indent{\it keywords:} backward doubly stochastic differential equations,
backward stochastic integral, comparison theorem
\end{abstract}

\section{Introduction}\label{sec:1}

Nonlinear backward stochastic differential equations (BSDEs in short)
have been independently introduced by Pardoux and Peng \cite{PP1} and Duffie
and Epstein \cite{DE}. Since then, BSDEs have been studied intensively. In
particular, many efforts have been made to relax the assumption on the
generator. For instance, Lepeltier and San Martin \cite{LS} have proved the
existence of a solution for the case when the generator is only continuous
with linear growth, and Jia \cite{J1,J2} studied
the existence of BSDEs with left-Lipschitz coefficients.
Another main reason is due to their enormous range of
applications in such diverse fields as mathematical finance (see Duffie
and Epstein \cite{DE} and Peng \cite{P2}), partial differential equations (see Peng \cite{P1}),
stochastic optimal control and stochastic game (see Hamadene and Lepeltier \cite{HL}),
nonlinear mathematical expectations (see Jiang and Chen \cite{JC} and Hu and Peng \cite{HP}), and so on.

A class of backward doubly stochastic differential equations (BDSDEs
in short) was introduced by Pardoux and Peng \cite{PP2} in 1994, in order
to provide a probabilistic interpretation for the solutions of a
class of semilinear stochastic partial differential equations (SPDEs
in short). They have proved the existence and uniqueness of
solutions for BDSDEs under uniformly Lipschitz conditions. Since
then, Shi et al. \cite{SGL} have relaxed the Lipschitz assumptions to
linear growth conditions. Bally and Matoussi \cite{BM} have given a
probabilistic interpretation of the solutions in Sobolev spaces for
semilinear parabolic SPDEs in terms of BDSDEs. Zhang and Zhao \cite{ZZ}
have proved the existence and uniqueness of solution for BDSDEs on
infinite horizons, and described the stationary solutions of SPDEs
by virtue of the solutions of BDSDEs on infinite horizons.
N'zi and Owo \cite{NO} have proved the existence
of a solution for one dimensional BDSDEs when the coefficient is
linear growth. Lin \cite{L} has also proved the existence of a solution
for one dimensional BDSDEs when the coefficient is bounded monotone.
Recently, Ren et al. \cite{RLH} and Hu and Ren \cite{HR} considered the
BDSDEs driven by Levy process with Lipschitz coefficient and
applications in SPDEs.

Unfortunately, most existence or uniqueness results of solution of BDSDEs
need the generator be at least continuous , which is somehow too strong in
some applications. Indeed, there are many SPDEs, in which the generator
may be discontinuous, and these PDEs have associated existence results of
solution (see Yoo \cite{Y} and Kim \cite{K}). Thus, a natural and interesting problem
is: can we establish the connections between SPDEs with discontinuous coefficient
and BDSDEs? Of course, the first step should be to obtain the existence and
uniqueness result of BDSDEs with discontinuous coefficient, next, to construct the
connections such as stochastic Feynman-Kac formula. Under which conditions do the
BDSDEs with discontinuous g have adapted solution?

Because of their important significance to SPDEs, it is necessary to
give intensive investigation to the theory of BDSDEs.
In this paper we shall study one-dimensional BDSDEs
\begin{eqnarray}\label{eq:1}
Y_t=\xi+\int_t^Tf(s,Y_s,Z_s)ds+\int_t^Tg(s,Y_s,Z_s)dB_s-
\int_t^TZ_sdW_s,\quad 0\leq t\leq T,
\end{eqnarray}
where $\xi\in L^2(\Omega,{\cal{F}}_T,P)$, $f:\Omega\times
[0,T]\times \mathbb{R}\times \mathbb{R}^{d}\rightarrow
\mathbb{R}$, and $g:\Omega\times [0,T]\times \mathbb{R}\times
\mathbb{R}^{d}\rightarrow \mathbb{R}^{l}$,  may be
discontinuous in $y$. Note
that the integral with respect to $\{B_t\}$ is a ``backward It\^o
integral" and the integral with respect to $\{W_t\}$ is a standard
forward It\^o integral. These two types of integrals are particular
cases of the It\^o-Sokorohod integral in Nualart and Pardoux \cite{NP}.
 In fact, we show that the one-dimensional
BDSDE associated with $(f,g,T,\xi)$ has at least a solution if $f$
and $g$ satisfy the following conditions:
\begin{enumerate}
\item[(H1)]
$f(t, \cdot, z)$ is left-continuous, and $f(t, y, \cdot)$ is
Lipschitz continuous, i.e., there exists a constant $K>0$, such
that $|f(t,y,z_1)-f(t,y,z_2)|\leq K|z_1-z_2|$, for all $t\in
[0,T]$, $y\in \mathbb{R}$, $z_1$, $z_2\in \mathbb{R}^d. $
\item[(H2)]
There exist two BDSDEs with generators $f_1$, $f_2$ respectively,
such that $f_1(t,y,z)\leq f(t,y,z)\leq f_2(t,y,z),$ for all $t\in
[0,T]$, $y\in \mathbb{R}$, $z\in \mathbb{R}^d$, and for given $T$
and $\xi$, the equations $(f_1,g,T,\xi)$ and $(f_2,g,T,\xi)$ have at
least one solution respectively, denoted by $(Y^i_t,Z^i_t)$,
$i=1,2$, where $Y^1_t\leq Y^2_t,$ for $t\in [0,T],$ a.s., a.e.
Moreover, the processes $f_i(t,Y^i_t,Z^i_t)$ are square integrable.
\item[(H3)]
$f(t,\cdot,z)$ satisfies left Lipschitz condition in $y,$ i.e.,
$f(t,y_1,z)-f(t,y_2,z)\geq -K(y_1-y_2),$ for all $t\in [0,T]$,
$y_1, y_2\in \mathbb{R}$, and $y_1\geq y_2, z\in \mathbb{R}^d. $
\item[(H4)]
There exist constants $c>0$ and $0<\alpha<1$ such that
 $|g(\omega,t,y_1,z_1)-g(\omega,t,y_2,z_2)|^2\leq c|y_1-y_2|^2+\alpha
|z_1-z_2|^2,$ for all $(\omega,t)\in \Omega\times [0,T],\
(y_1,z_1)\in \mathbb{R}\times \mathbb{R}^d, \ (y_2,z_2)\in
\mathbb{R}\times \mathbb{R}^d.$
\end{enumerate}

It should be noted that our conditions of this paper, without explicit growth constraint,
is different from N'zi and Owo \cite{NO};  and without the monotone and bounded constraint,
is different from Lin \cite{L}.

This paper is organized as follows. In Section \ref{sec:2} we formulate the
problem accurately and give some preliminary results.
Section \ref{sec:3} is devoted to the proof of the existence of solutions of BDSDEs.
Finally, in Section \ref{sec:4} the comparison theorem is obtained.

\section{Preliminaries}\label{sec:2}

Let $(\Omega,{\cal{F}},P)$ be a probability space, and $T>0$ be an
arbitrarily fixed constant throughout this paper. Let $\{W_t;0\leq
t\leq T\}$ and $\{B_t;0\leq t\leq T\}$ be two mutually independent
standard Brownian Motions with values in $\mathbb{R}^d$ and
$\mathbb{R}^l$, respectively, defined on $(\Omega,{\cal{F}},P)$. Let
$\cal{N}$ denote the class of $P$-null sets of $\cal{F}$. For each
$t\in [0,T]$, we define ${\cal{F}}_t={\cal{F}}_t^W\vee
{\cal{F}}_{t,T}^B$, where for any process $\{\eta_t\}$,
${\cal{F}}_{s,t}^{\eta}=\sigma\{\eta_r-\eta_s;s\leq r\leq t\}\vee
{\cal{N}}$, ${\cal{F}}_t^{\eta}={\cal{F}}_{0,t}^{\eta}$. Note that
$\{{\cal{F}}_t;t\in [0,T]\}$ is neither increasing nor decreasing,
so it does not constitute a common filtration. We introduce the
following notations:
\begin{eqnarray*}
S^2\left([0,T];\mathbb{R}\right)&=&\{v_t,0\leq t\leq T,\
 \mbox {is a}\ \mathbb{R}\mbox{-valued},\ {\cal F}_t \mbox{-measurable continuous }\\
 &&  \mbox{process such that}\ E(\sup_{0\leq t\leq T}|v_{t}|^{2})<\infty\},\\
M^2(0,T;\mathbb{R}^n)&=&\{v_t,0\leq t\leq T,\
 \mbox {is a}\ \mathbb{R}^n\mbox{-valued},\ {\cal F}_t\mbox{-measurable process}\\
 &&  \mbox{such that}\ E\int_{0}^T|v_{t}|^{2}dt<\infty\}.
\end{eqnarray*}
We use the usual inner product $\langle \cdot ,\cdot \rangle $ and
Euclidean norm $ | \cdot | $ in $\mathbb{R}$, $\mathbb{R}^{l}$ and
$\mathbb{R}^{d}.$ All the equalities and inequalities mentioned in
this paper are in the sense of $dt\times dP$ almost
surely on $\left[0,T\right] \times \Omega$.

\begin{definition}\label{def:2.1}
A pair of processes $(y,z):\Omega \times
[0,T]\rightarrow \mathbb{R}\times \mathbb{R}^{d}$ is called a solution of
BDSDE (\ref{eq:1}) if $(y,z)\in S^2([0,T];\mathbb{R})\times
M^2(0,T;\mathbb{R}^{d})$ and satisfies BDSDE (\ref{eq:1}).
\end{definition}
Also we need one lemma, which is a special case of
comparison theorem in Shi et al. \cite{SGL} .

\begin{lemma}\label{lem:2.2}
Let $f_1(s,y,z)=ly+m|z|$, $f_2(s,y,z)=l|y|+m|z|$,
where constants $l, m \in \mathbb{R}$, and positive process $\phi
\in M^2(0,T;\mathbb{R})$, furthermore, $(y_t^i,z_t^i)_{t\in[0,T]}$
 $(i=1,2)$ are the solution to the following equations:
\begin{eqnarray*}
y^i_t=\xi+\int_t^T(f_i(s,y^i_s,z^i_s)+\phi_s)ds
+\int_t^Tg(s,y^i_s,z^i_s)dB_s-\int_t^Tz^i_sdW_s,
\ i=1,2.
\end{eqnarray*}
If $\xi\in L^2(\Omega,{\cal F}_T,P)$ and $\xi\geq 0$, then
$y_t^i\geq 0$, P-a.s., $t\in[0,T], i=1,2.$
\end{lemma}
\begin{remark}\label{rmk:2.3}
The assumptions (H1) and (H3) imply
\begin{eqnarray*}f(t,y_1,z_1)-f(t,y_2,z_2)\geq -K(y_1-y_2)-K|z_1-z_2|,
\hspace{1mm}y_1, y_2\in\mathbb{R}, z_1, z_2\in\mathbb{R}^d.
\end{eqnarray*}
\end{remark}

\section{Existence of Solutions of BDSDEs}\label{sec:3}

In this section, we will state and prove the existence of solutions of BDSDEs.
\begin{theorem}\label{thm:3.1}
Under the assumptions (H1)-(H4), then, if $\xi \in
L^2(\Omega,{\cal{F}}_T,P)$, BDSDE (\ref{eq:1}) has a
solution $(Y_t,Z_t) \in S^2([0,T]$; $\mathbb{R})\times
M^2(0,T;\mathbb{R}^d).$
\end{theorem}
At first, we denote that $(Y^j_t,Z^j_t)$ are the solutions of
$(f_j, g, T, \xi),$ where $j=1,2,$ that is
\begin{eqnarray}\label{eq:2}
Y^j_t=\xi+\int_t^Tf_j(s,Y^j_s,Z^j_s)ds+\int_t^Tg(s,Y^j_s,Z^j_s)dB_s-\int_t^TZ^j_sdW_s,
\end{eqnarray}
where $f_j$ satisfies (H2) and $f_j(t,Y^j_t,Z^j_t)\in
M^2(0,T;\mathbb{R}).$ Now we construct a sequence of BDSDEs as
follows:
\begin{eqnarray}\label{eq:3}
\nonumber y^{i}_t&=&\xi+\int_t^T\left(f(s,y^{i-1}_s,z^{i-1}_s)
-K(y^{i}_s-y^{i-1}_s)-K|z^{i}_s-z^{i-1}_s|\right)ds\\
&&+\int_t^Tg(s,y^{i}_s,z^{i}_s)dB_s-\int_t^Tz^{i}_sdW_s,
\end{eqnarray}
where $i=1,2,\cdots$, and $(y^{0}_t,z^{0}_t)=(Y^1_t,Z^1_t).$
By Theorem 1.1 in Pardoux and Peng \cite{PP2},
BDSDEs (\ref{eq:3}) $(i=1,2,\cdots)$ have a unique
adapted solution respectively if $f(t,y^{i-1}_t,z^{i-1}_t)\in
M^2(0,T;\mathbb{R}).$ For these equations, we have:

\begin{lemma}\label{lem:3.2}
Under the assumptions (H1)-(H4), for any positive integer
$i$, then, if $\xi \in L^2(\Omega,{\cal{F}}_T,P)$,
BDSDE (\ref{eq:3}) has a unique adapted solution
$(y^{i}_t,z^{i}_t)\in S^2([0,T]$; $\mathbb{R})\times
M^2(0,T;\mathbb{R}^{d})$ and $Y^1_t\leq y^{i}_t\leq y^{i+1}_t\leq
Y^2_t$, P-a.s., $\forall t\in [0,T].$
\end{lemma}

\noindent{\bf Proof.}\quad
Firstly, for $i=1,$ by $Y^2_t\geq Y^1_t$ and (H2), it follows that
\begin{eqnarray}\label{eq:4}
f_2(t,Y^2_t,Z^2_t)-f(t,Y^1_t,Z^1_t)
&\geq&
f(t,Y^2_t,Z^2_t)-f(t,Y^1_t,Z^1_t)\\
&\geq&
-K(Y^{2}_t-Y^{1}_t)-K|Z^{2}_t-Z^{1}_t|.
\end{eqnarray}
Thus
\begin{eqnarray*}
f_2(t,Y^2_t,Z^2_t)+K(Y^{2}_t-Y^{1}_t)+K|Z^{2}_t-Z^{1}_t| \geq
f(t,Y^1_t,Z^1_t) \geq f_1(t,Y^1_t,Z^1_t),
\end{eqnarray*}
this implies that $f(t,Y^1_t,Z^1_t) \in M^2(0,T;\mathbb{R})$. By
Theorem 1.1 in Pardoux and Peng \cite{PP2}, it follows that BDSDE (\ref{eq:3})
has a unique adapted solution $(y^{1}_t,z^{1}_t).$

Now, by (\ref{eq:3}) and (\ref{eq:2}) when $i=1$ and $j=1$, we have
\begin{eqnarray}
\nonumber y^1_t-Y^1_t&=&\int_t^T(-K(y^1_s-Y^1_s)-K|z^1_s-Z^1_s|+\Delta_s^1)ds\\
\nonumber &&+\int_t^T(g(s,y^1_s,z^1_s)-g(s,Y_s^1,Z_s^1))dB_s-\int_t^T(z_s^1-Z_s^1)dW_s,
\end{eqnarray}
where
$\Delta_s^1:=f(s,Y_s^{1},Z_s^{1})-f_1(s,Y_s^{1},Z_s^{1}) \geq 0$ and $\Delta_s^1
\in M^2(0,T;\mathbb{R}).$ Applying Lemma \ref{lem:2.2}, we know
$y^1_t\geq Y^1_t,\hspace{0.1cm}
\mbox{P-a.s.,}\hspace{0.1cm}\forall t\in [0,T]$.

Again we consider (\ref{eq:3}) and (\ref{eq:2}) when $i=1$ and $j=2$, we have
\begin{eqnarray}
\nonumber Y^2_t-y^1_t&=&\int_t^T(-K(Y^2_s-y^1_s)-K|Z^2_s-z^1_s|+\Delta_s^2)ds\\
\nonumber &&+\int_t^T(g(s,Y^2_s,Z^2_s)-g(s,y_s^1,y_s^1))dB_s-\int_t^T(Z_s^2-z_s^1)dW_s,
\end{eqnarray}
where
$\Delta_s^2:=f_2(s,Y_s^{2},Z_s^{2})-f(s,Y_s^{1},Z_s^{1})+K(Y^2_s-Y^1_s)
+K|Z^2_s-Z^1_s| \geq 0$ and $\Delta_s^2
\in M^2(0,T;\mathbb{R}).$
Then from Lemma \ref{lem:2.2}, we have
\begin{eqnarray*}
Y^1_t\leq y^{1}_t\leq Y^{2}_t,\ \mbox{P-a.s.,}\ \forall t \in [0,T].
\end{eqnarray*}

Similarly, for $i=2,$ since $Y^1_t\leq y^{1}_t\leq Y^{2}_t$ and
(H2), it follows that
\begin{eqnarray*}
f_2(t,Y^2_t,Z^2_t)-f(t,y^{1}_t,z^{1}_t) &\geq&
f(t,Y^2_t,Z^2_t)-f(t,y^{1}_t,z^{1}_t)\\
&\geq&
-K(Y^{2}_t-y^{1}_t)-K|Z^{2}_t-z^{1}_t|.
\end{eqnarray*}
Thus
\begin{eqnarray*}
f_2(t,Y^2_t,Z^2_t)+K(Y^{2}_t-y^{1}_t)+K|Z^{2}_t-z^{1}_t| \geq
f(t,y^{1}_t,z^{1}_t).
\end{eqnarray*}
But
\begin{eqnarray*}
f(t,y^{1}_t,z^{1}_t) \geq
f_1(t,Y^1_t,Z^1_t)-K(y^{1}_t-Y^{1}_t)-K|z^{1}_t-Z^{1}_t|,
\end{eqnarray*}
this implies that $f(t,y^{1}_t,z^{1}_t) \in M^2(0,T;\mathbb{R})$, and
BDSDE (\ref{eq:3}) has a unique adapted solution
$(y^{2}_t,z^{2}_t).$ Using the similar method, we have
\begin{eqnarray*}
y^{1}_t\leq y^{2}_t\leq Y^{2}_t,\ \mbox{P-a.s.,}\ \forall t \in
[0,T].
\end{eqnarray*}

Finally, for $i>2$, we assume that $Y^1_t\leq y^{i-1}_t\leq
y^{i}_t\leq Y^2_t$ and $f(t,y^{i-1}_t,z^{i-1}_t) \in M^2(0,T;R),$ we
consider BDSDE (\ref{eq:3}) for $i+1,$ which can be written as
\begin{eqnarray}\label{eq:5}
\nonumber y^{i+1}_t&=&\xi+\int_t^T\left(f(s,y^{i}_s,z^{i}_s)
-K(y^{i+1}_s-y^{i}_s)-K|z^{i+1}_s-z^{i}_s|\right)ds\\
&&+\int_t^Tg(s,y^{i+1}_s,z^{i+1}_s)dB_s-\int_t^Tz^{i+1}_sdW_s.
\end{eqnarray}
Then by the similar argument as the case $i=2,$ we have
\begin{eqnarray*}
&&f_2(t,Y^2_t,Z^2_t)+K(Y^{2}_t-y^{i}_t)+K|Z^{2}_t-z^{i}_t|\\
&\geq& f(t,y^{i}_t,z^{i}_t)\geq
f_1(t,Y^1_t,Z^1_t)-K(y^{i}_t-Y^{1}_t)-K|z^{i}_t-Z^{1}_t|,
\end{eqnarray*}
this implies that $f(t,y^{i}_t,z^{i}_t) \in M^2(0,T;\mathbb{R})$, and
BDSDE (\ref{eq:5}) has a unique adapted solution
$(y^{i+1}_t,z^{i+1}_t).$ By  Lemma \ref{lem:2.2} again, we have
\begin{eqnarray*}
Y^{1}_t\leq y^{i}_t\leq y^{i+1}_t\leq Y^{2}_t,\ \mbox{P-a.s.,}\
\forall t \in [0,T].
\end{eqnarray*}
\quad$\Box$

\begin{lemma}\label{lem:3.3}
There exists a constant $A>0$, such that
\begin{eqnarray*}
\sup\limits_iE\left[\sup\limits_{0\leq t\leq
T}|y_t^i|^2+\int_0^T|z_t^i|^2dt\right]<A.
\end{eqnarray*}
\end{lemma}
{\bf Proof.}\quad From Lemma \ref{lem:3.2}, we have
\begin{eqnarray*}
\sup\limits_i\left[E\sup\limits_{0\leq t\leq T}|y_t^i|^2\right]
\leq E\left[\sup\limits_{0\leq t\leq T}|Y_t^1|^2\right]+
E\left[\sup\limits_{0\leq t\leq T}|Y_t^2|^2\right]<\infty.
\end{eqnarray*}
By the similar argument as (\ref{eq:4}), we can deduce
\begin{eqnarray*}
&&f_2(t,Y^2_t,Z^2_t)+K(Y^{2}_t-y^{i}_t)+K|Z^{2}_t-z^{i}_t|\\
&\geq&f(t,y^i_t,z^i_t)\geq
f_1(t,Y^1_t,Z^1_t)-K(y^{i}_t-Y^{1}_t)-K|z^{i}_t-Z^{1}_t|.
\end{eqnarray*}
Then, we have
\begin{eqnarray}\label{eq:6}
\nonumber
&&|f(t,y^i_t,z^i_t)-K(y^{i+1}_t-y^{i}_t)-K|z^{i+1}_t-z^{i}_t||\\
\nonumber
&\leq& |f(t,y^i_t,z^i_t)|+|K(y^{i+1}_t-y^{i}_t)|+K|z^{i+1}_t-z^{i}_t|\\
\nonumber
&\leq&|f_1(t,Y^1_t,Z^1_t)-K(y^{i}_t-Y^{1}_t)-K|z^{i}_t-Z^{1}_t||+K|y^{i+1}_t-y^{i}_t|\\
\nonumber
&&+|f_2(t,Y^2_t,Z^2_t)-K(Y^{2}_t-y^{i}_t)-K|Z^{2}_t-z^{i}_t||+K|z^{i+1}_t-z^{i}_t|\\
\nonumber
&\leq& \sum\limits_{j=1}^2\left[|f_j(t,Y^i_t,Z^i_t)|+K(|Y^{j}_t|+|Z^{j}_t|)\right]\\
&&+3K(|y^{i}_t|+|z^{i}_t|)+K(|y^{i+1}_t|+|z^{i+1}_t|).
\end{eqnarray}
Applying It\^o's formula to $|y^{i+1}_t|^2$ for $t\in [0,T]$, we
deduce
\begin{eqnarray*}
&&E\int_0^T|z^{i+1}_t|^2dt\\
&\leq& E|\xi|^2+2E\int_0^Ty^{i+1}_t\cdot
(f(t,y^i_t,z^i_t)-K(y^{i+1}_t-y^{i}_t)-K|z^{i+1}_t-z^{i}_t|)dt\\
&&+E\int_0^T|g(t,y^{i+1}_t,z^{i+1}_t)|^2dt\\
&\leq& E|\xi|^2+2E\int_0^T|y^{i+1}_t|\cdot
|f(t,y^i_t,z^i_t)-K(y^{i+1}_t-y^{i}_t)-K|z^{i+1}_t-z^{i}_t||dt\\
&&+E\int_0^T|g(t,y^{i+1}_t,z^{i+1}_t)|^2dt\\
&\leq& C_1+\dfrac{1-\alpha'}{4}E\int_0^T(|z^{i}_t|^2+|z^{i+1}_t|^2)dt
+E\int_0^T|g(t,y^{i+1}_t,z^{i+1}_t)|^2dt,
\end{eqnarray*}
\noindent with some constant $C_1>0$. Hereafter, $\forall n \geq 1,\
C_n$ will be some positive real constant. Then
\begin{eqnarray*}
E\int_0^T|z^{i+1}_t|^2dt \leq
C_2+\dfrac{1-\alpha'}{4}E\int_0^T(|z^{i}_t|^2+|z^{i+1}_t|^2)dt +\alpha'E\int_0^T|z^{i+1}_t|^2dt.
\end{eqnarray*}
That is
\begin{eqnarray*}
E\int_0^T|z^{i+1}_t|^2dt\leq
C_3+\dfrac{1}{3}E\int_0^T|z^{i}_t|^2dt.
\end{eqnarray*}
This implies that $\sup\limits_iE\int_0^T|z^{i+1}_t|^2dt <\infty$,
which yields that the quantities
$\psi^{i+1}(t$, $y^{i+1}_t,z^{i+1}_t)=f(t,y^i_t,z^i_t)$
$-K(y^{i+1}_t-y^{i}_t)-K|z^{i+1}_t-z^{i}_t|$
are uniformly bounded in $M^2(0,T;\mathbb{R}).$ Set
$C_0=\sup\limits_iE\int_0^T|\psi^{i+1}(t,y^{i+1}_t$,
$z^{i+1}_t)|^2dt.$ \quad$\Box$

\begin{lemma}\label{lem:3.4}
There exist processes $(y_t,z_t)\in S^2([0,T];\mathbb{R})\times
M^2(0,T;\mathbb{R}^{d})$  such that as $n\rightarrow \infty$
$$E\left[\sup\limits_{0\leq t\leq T}|y_t^n-y_t|^2
+\int_0^T|z_t^n-z_t|^2dt\right]\rightarrow 0.$$
\end{lemma}
{\bf Proof.}\quad By Lemma \ref{lem:3.2}, it follows that there
exists a process $y_t$ such that
$
y_t^n \nearrow y_t,\ \mbox{P-a.s.},\ \forall t\in [0,T],
$
as $n \rightarrow \infty$, and $E[\sup\limits_{0\leq t\leq
T}|y_t|^2]<\infty.$ By the dominated convergence theorem, we get as
$n \rightarrow \infty,$
\begin{eqnarray}\label{eq:7}
E\int_0^T|y^n_t-y_t|^2dt\rightarrow 0.
\end{eqnarray}
Coming back to (\ref{eq:6}) and by Lemma \ref{lem:3.2}, we can deduce
\begin{eqnarray}\label{eq:8}
\sup\limits_{n}\left[E\int_0^T|\psi^n(t,y_t^n,z_t^n)|^2dt\right]<\infty.
\end{eqnarray}
Applying It\^o's formula to $|y^n_t-y^m_t|^2$ for $t\in [0,T]$,
taking expectation in both sides, we have
\begin{eqnarray*}
& &E(|y^n_0-y^m_0|^2)+E\int_0^T|z^n_t-z^m_t|^2dt \\
&=&2E\int_0^T(y^n_t-y^m_t)(\psi^{n}(t,y^{n}_t,z^{n}_t)-\psi^{m}(t,y^{m}_t,z^{m}_t))dt\\
&&+E\int_0^T|g(t,y^n_t,z^n_t)-g(t,y^m_t,z^m_t)|^2dt\\
&\leq& 4C_0(E\int_0^T |y^n_t-y^m_t|^2dt)^{1/2}
+E\int_0^T (c|y^n_t-y^m_t|^2+\alpha |z^n_t-z^m_t|^2 )dt.
\end{eqnarray*}
Then
\begin{eqnarray*}
&&E\int_0^T|z^n_t-z^m_t|^2dt\\
&\leq& \dfrac{1}{1-\alpha}\{4C_0(E\int_0^T |y^n_t-y^m_t|^2dt)^{1/2}
+cE\int_0^T|y^n_t-y^m_t|^2dt\}.
\end{eqnarray*}
Thus $\{z_t^n\}$ is a Cauchy sequence in $M^2(0,T;\mathbb{R}^d)$,
therefore $\{z^{n}_t\}_{n=1}^\infty$ converges in
$M^2(0,T;\mathbb{R}^d)$, to a limit $z_t$, we have
\begin{eqnarray}\label{eq:9}
E\int_0^T|z^n_t-z_t|^2dt\rightarrow 0.
\end{eqnarray}
Applying It\^o's formula to $|y^n_s-y^m_s|^2$ for $s\in [t,T]$,
we have
\begin{eqnarray}
\nonumber &&|y^n_t-y^m_t|^2+\int_t^T|z^n_s-z^m_s|^2ds \\
\nonumber &=&
2\int_t^T(y^n_s-y^m_s)(\psi^{n}(s,y^{n}_s,z^{n}_s)-\psi^{m}(s,y^{m}_s,z^{m}_s))ds\\
\nonumber &&+\int_t^T|g(s,y^n_s,z^n_s)-g(s,y^m_s,z^m_s)|^2ds
-2\int_t^T(y^n_s-y^m_s)\cdot(z^n_s-z^m_s)dW_s\\
\nonumber&&
+2\int_t^T(y^n_s-y^m_s)\cdot (g(s,y^n_s,z^n_s)-g(s,y^m_s,z^m_s))dB_s.
\end{eqnarray}
\noindent Taking supremum and expectation, by Young's inequality, we
get
\begin{eqnarray}\label{eq:10}
\nonumber
&&E\left[\sup\limits_{0\leq t\leq T}|y^n_t-y^m_t|^2\right]\\
\nonumber
&\leq& 2\left[E\int_0^T(y^n_s-y^m_s)^2ds\right]^{1/2}
\left[E\int_0^T(\psi^{n}(s,y^{n}_s,z^{n}_s)-\psi^{m}(s,y^{m}_s,z^{m}_s))^2ds\right]^{1/2}\\
\nonumber
&&+E\int_0^T|g(s,y^n_s,z^n_s)-g(s,y^m_s,z^m_s)|^2ds\\
\nonumber
&&+2E\sup_{0\leq t\leq T}\left|\int_t^T(y^n_s-y^m_s)\cdot(z^n_s-z^m_s)dW_s\right|\\
& & +2E\sup_{0\leq t\leq T}\left|\int_t^T(y^n_s-y^m_s)\cdot
(g(s,y^n_s,z^n_s)-g(s,y^m_s,z^m_s))dB_s\right|.
\end{eqnarray}
By Burkholder-Davis-Gundy's inequality, we deduce
\begin{eqnarray}\label{eq:11}
\nonumber
&& E(\sup_{0\leq t\leq T}|\int_t^T(y^n_s-y^m_s)\cdot
(g(s,y^n_s,z^n_s)-g(s,y^m_s,z^m_s))dB_s|)\\
\nonumber
& \leq &  k E(\int_0^T|y^n_s-y^m_s|^2\cdot
|g(s,y^n_s,z^n_s)-g(s,y^m_s,z^m_s))|^2ds)^{1/2}\\
\nonumber
& \leq & k E((\sup_{0\leq t\leq T }|y^n_t-y^m_t|^2)^{1/2}
( \int_0^T |g(s,y^n_s,z^n_s)-g(s,y^m_s,z^m_s)|^2 ds)^{1/2})\\
\nonumber
& \leq &2k^2cE \int_0^T|y^n_s-y^m_s|^2ds
+\displaystyle\frac{1}{8}E(\sup\limits_{0\leq t\leq T}
|y^n_t-y^m_t|^2)\\
&&+2k^2\alpha E\int_0^T|z^n_s-z^m_s|^2ds.
\end{eqnarray}
In the same way, we have
\begin{eqnarray}\label{eq:12}
\nonumber
&&E(\sup_{0\leq t\leq T}|\int_t^T(y^n_s-y^m_s)\cdot(z^n_s-z^m_s)dW_s|)\\
&&\leq\displaystyle\frac{1}{8}E(\sup\limits_{0\leq t\leq T}
|y^n_t-y^m_t|^2)+2k^2\int_0^T|z^n_s-z^m_s|^2ds.
\end{eqnarray}
From (\ref{eq:11}), (\ref{eq:12}) and (\ref{eq:10}), it follows that
\begin{eqnarray*}
&&E\left[\sup\limits_{0\leq t\leq T}|y^n_t-y^m_t|^2\right]\\
&\leq&4C_0\left[E\int_0^T(y^n_s-y^m_s)^2ds\right]^{1/2}
+2c(4k^2+1)E\int_0^T|y^n_s-y^m_s)|^2ds\\
&& +2(\alpha+4k^2(\alpha+1))\int_0^T|z^n_s-z^m_s|^2ds.
\end{eqnarray*}
Then from (\ref{eq:7}) and (\ref{eq:9}), we can deduce
\begin{eqnarray*}
E\left[\sup\limits_{0\leq t\leq T}|y^m_t-y^n_t|^2\right]\rightarrow 0
\end{eqnarray*}
 as $n, m\rightarrow \infty.$ Obviously, the process $y_t$
 belongs to $S^2([0,T];\mathbb{R}).$

By (H1), (\ref{eq:7}) and (\ref{eq:9}), it follows that there
exists a subsequence (we still denote by $n$) such that as $n
\rightarrow \infty,$
\begin{eqnarray}\label{eq:13}
\psi^n(t,y^n_t,z^n_t)-f(t,y_t,z_t) \rightarrow 0,\ dt\times dP\mbox{-a.s.}\
\end{eqnarray}
\quad$\Box$

Now we are in the position to give the proof of Theorem \ref{thm:3.1}.

\noindent\textbf{Proof of Theorem \ref{thm:3.1}}\quad By
(\ref{eq:9}) and the continuity of the stochastic integral, we get
\begin{eqnarray*}
&&\sup\limits_{0\leq t\leq T} \left|\int_t^Tz^n_sdW_s-
\int_t^Tz_sdW_s\right|\rightarrow 0
\hspace{0.3cm} \mbox{in probability,}\\
&&\sup\limits_{0\leq t\leq T} \left|\int_t^Tg(s,y_s^n,z_s^n)dB_s
-\int_t^Tg(s,y_s,z_s)dB_s\right|\rightarrow
0\hspace{0.3cm} \mbox{in probability.}
\end{eqnarray*}
So there exist a subsequence (we still denote by $\{n\}$) such that
as the convergence is $P$-almost surely.

Since (H1), (H3) and $Y_t^1\leq y_t\leq Y_t^2,$ we can deduce from
the similar argument as (\ref{eq:4}) that $f(t,y_t,z_t)\in
M^2([0,T];\mathbb{R})$. In view of (\ref{eq:9}) and
(\ref{eq:13}), then by the dominated convergence theorem, passing
to a subsequence (we still denote by $\{n\}$), we have
\begin{eqnarray*}
\int_0^T|\psi^n(t,y^n_t,z^n_t)-f(t,y_t,z_t)|dt\rightarrow 0,\quad
\mbox{P-a.s.}.
\end{eqnarray*}
Hence, passing to the limit, as $ i \rightarrow \infty$ on both
sides of (\ref{eq:3}), we can get
\begin{eqnarray*}
y_t=\xi+\int_t^Tf(s,y_s,z_s)ds+\int_t^Tg(s,y_s,z_s)dB_s-\int_t^Tz_sdW_s.
\end{eqnarray*}
It is obvious that $(y_t,z_t)$ is a solution of
BDSDE (\ref{eq:1}) under (H1)-(H4). \quad$\Box$

\begin{remark}\label{rmk:3.5}
Although the solution $(y_t,z_t)$ we get
in Theorem \ref{thm:3.1} is constructed by approximating from below,
we cannot get that the constructed solution is the minimal
solution of BDSDE (\ref{eq:1}), in the sense that for any other
solution $(Y_t,Z_t)$ of BDSDE (\ref{eq:1}), we have $y_t\leq Y_t$
under (H1)-(H4). This is because we cannot compare the
generators of BDSDE (\ref{eq:1}) and BDSDE (\ref{eq:3}) such that
Lemma \ref{lem:2.2} cannot be used to compare the solutions
of BDSDE (\ref{eq:1}) and BDSDE (\ref{eq:3}).
\end{remark}

In order to get the minimal solution of BDSDE (\ref{eq:1}),
in the following of this paper, we replace (H2) by
\begin{enumerate}
\item[(H5)]
$|f(t,y,0)|\leq |f(t,0,0)|+K|y|,\\ \quad\forall ( t, y, z)\in
[0,T]\times \mathbb{R}\times \mathbb{R}^d$, and $
E\int_0^T|f(t,0,0)|^2{\rm dt}<\infty. $
\end{enumerate}

\begin{remark}\label{rmk:3.6}
\begin{enumerate}
\item[(i)]
The assumptions (H1) and (H5) imply
\begin{eqnarray*}
|f(t,y,z)|\leq K|y|+K|z|+|f(t,0,0)|,
\hspace{1mm}y \in\mathbb{R}, z\in\mathbb{R}^d, t\in
[0,T].
\end{eqnarray*}
\item[(ii)]
It is obvious that (H5) is a special case of (H2).
\end{enumerate}
\end{remark}

We construct a sequence of BDSDEs as follows:
\begin{eqnarray}\label{eq:14}
\nonumber\underline{Y}_t^0&=&\xi+\int_t^T(-K|\underline{Y}_s^0|
-K|\underline{Z}_s^0|-|f(s,0,0)|)ds\\
&&+\int_t^Tg(s,\underline{Y}_s^0,\underline{Z}_s^0)dB_s
-\int_t^T\underline{Z}_s^0dW_s,
\end{eqnarray}
\begin{eqnarray}\label{eq:15}
\nonumber\underline{Y}_t^{i+1}&=&\xi+\int_t^T(f(s,\underline{Y}_s^{i},
\underline{Z}_s^{i})-K(\underline{Y}_s^{i+1}-\underline{Y}_s^{i})-
K|\underline{Z}_s^{i+1}-\underline{Z}_s^{i}|)ds\\
&&+\int_t^Tg(s,\underline{Y}_s^{i+1},\underline{Z}_s^{i+1})dB_s
-\int_t^T\underline{Z}_s^{i+1}dW_s,
\end{eqnarray}
where $i=0,1,2,\cdots.$ Besides the above equations, we also need
the other BDSDE
\begin{eqnarray}\label{eq:16}
\nonumber
\overline{Y}_t^0&=&\xi+\int_t^T(K|\overline{Y}_s^0|
+K|\overline{Z}_s^0|+|f(s,0,0)|)ds\\
&&+\int_t^Tg(s,\overline{Y}_s^0,\overline{Z}_s^0)dB_s
-\int_t^T\overline{Z}_s^0dW_s.
\end{eqnarray}
By Theorem 1.1 in Pardoux and Peng \cite{PP2}, BDSDEs (\ref{eq:14}), (\ref{eq:15})
$(i=0,1,\cdots)$ and (\ref{eq:16}) have a unique adapted solution
respectively. For these solutions mentioned above, with the technique similarly
to Lemma \ref{lem:3.2}, we can obtain the following properties:

\begin{lemma}\label{lem:3.7}
Under the assumptions (H1) and (H3)-(H5), the following
properties hold true:
\begin{enumerate}
\item[(i)]
For any positive integer $i$,
$\underline{Y}^{i+1}_t\geq\underline{Y}^{i}_t\geq\underline{Y}^0_t $,
P-a.s., $\forall t\in [0,T].$
\item[(ii)]
For any positive integer $i$, $ \underline{Y}^{i}_t\leq
\overline{Y}^0_t$, P-a.s., $\forall t\in [0,T].$
\end{enumerate}
\end{lemma}

Lemma \ref{lem:3.7} implies that the sequence of solutions of BDSDEs
(\ref{eq:14}) and (\ref{eq:15}) is increasing and have upper bound
by the solution of (\ref{eq:16}), that is
\begin{eqnarray}\label{eq:17}
\underline{Y}^0_t\leq \underline{Y}^{i}_t\leq \underline{Y}^{i+1}_t\leq
\overline{Y}^0_t,\ \mbox{P-a.s.},\ \forall t\in [0,T],\
i=1,2,\cdots.
\end{eqnarray}
Furthermore, we get the existence of the minimal solution of BDSDEs.

\begin{theorem}\label{thm:3.8}
Under the assumptions (H1) and (H3)-(H5), then, if $\xi \in
L^2(\Omega$, ${\cal{F}}_T,P)$, BDSDE (\ref{eq:1}) has a solution
$(Y_t,Z_t) \in S^2([0,T];\mathbb{R})\times M^2(0,T;\mathbb{R}^d).$
Moreover, there is a minimal solution $(\underline{Y}_t, \underline{Z}_t)$
of BDSDE (\ref{eq:1}) in the sense that,
for any other solution $(Y_t,Z_t)_{t\in[0,T]}$ of BDSDE (\ref{eq:1}),
 we have $Y_t \geq \underline{Y}_t$, P-a.s., $\forall
t\in [0,T].$
\end{theorem}
{\bf Proof.}\ Similarly to the arguments in Theorem \ref{thm:3.1},
the existence of a solution of (\ref{eq:1}) can be obtained easily.
We want to prove the existence of a minimal solution of (\ref{eq:1}).
Let $(Y_t,Z_t)_{t\in[0,T]}$ be a solution of BDSDEs
(\ref{eq:1}).
For $i=0$, we have
\begin{eqnarray*}
Y_t-\underline{Y}^0_t&=&\int_t^T(f(s,Y_s,Z_s)
+K|\underline{Y}_s^0|+K|\underline{Z}_s^0|
+|f(s,0,0)|)ds\\
&&+\int_t^T(g(s,Y_s,Z_s)-g(s,\underline{Y}_s^0,\underline{Z}_s^0))dB_s
-\int_t^T(Z_s-\underline{Z}_s^0)dW_s\\
&=&\int_t^T(-K|Y_s-\underline{Y}_s^{0}|
-K|Z_s-\underline{Z}_s^{0}|+\Psi_s^0)ds\\
&&+\int_t^T(g(s,Y_s,Z_s)-g(s,\underline{Y}_s^0,\underline{Z}_s^0))dB_s
-\int_t^T(Z_s-\underline{Z}_s^0)dW_s
\end{eqnarray*}
where \begin{eqnarray*} \Psi_s^0&:=&K|Y_s-\underline{Y}_s^{0}|
+K|Z_s-\underline{Z}_s^{0}|+K|\underline{Y}_s^0|+K|\underline{Z}_s^0|
+|f(s,0,0)|-f(s,Y_s,Z_s)\\
&\geq& K|Y_s|+K|Z_s| +|f(s,0,0)|-f(s,Y_s,Z_s)\geq 0,
\end{eqnarray*}
From Lemma \ref{lem:2.2}, we have $Y_t\geq \underline{Y}^{0}_t$,
P-a.s., $\forall t\in [0,T]$.

We assume that $Y_t\geq \underline{Y}^{i}_t$, we will prove that
$Y_t\geq \underline{Y}^{i+1}_t$. We have
\begin{eqnarray*}
Y_t-\underline{Y}^{i+1}_t&=&\int_t^T(-K|Y_s-\underline{Y}_s^{i+1}|
-K|Z_s-\underline{Z}_s^{i+1}|+\Psi_s^{i+1})ds\\
&&+\int_t^T(g(s,Y_s,Z_s)-g(s,\underline{Y}_s^{i+1},\underline{Z}_s^{i+1}))dB_s
-\int_t^T(Z_s-\underline{Z}_s^{i+1})dW_s
\end{eqnarray*}
where
\begin{eqnarray*}
\Psi_s^{i+1}&:=&f(s,Y_s,Z_s)-f(s,\underline{Y}_s^{i},\underline{Z}_s^{i})
+K(\underline{Y}_s^{i+1}-\underline{Y}_s^{i})+K|\underline{Z}_s^{i+1}-\underline{Z}_s^{i}|\\
&&+K(\underline{Y}_s^{i+1}-Y_s)+K|\underline{Z}_s^{i+1}-Z_s|\\
&\geq &f(s,Y_s,Z_s)-f(s,\underline{Y}_s^{i},\underline{Z}_s^{i})
+K(Y_s-\underline{Y}_s^{i})+K|Z_s-\underline{Z}_s^{i}|\geq 0,
\end{eqnarray*}
then, $Y_t\geq \underline{Y}^{i+1}_t$, P-a.s., $\forall t\in [0,T]$.
This implies $Y_t \geq \underline{Y}_t$, P-a.s., $\forall t\in
[0,T].$ \quad$\Box$

In order to get the upper bound of solution of BDSDE (\ref{eq:1}),
besides (\ref{eq:14}), we also need the following BDSDE:
\begin{eqnarray}\label{eq:18}
\overline{Y}_t^{i+1}&=&\xi+\int_t^T(f(s,\overline{Y}_s^{i},
\overline{Z}_s^{i})-K(\overline{Y}_s^{i+1}-\overline{Y}_s^{i})+
K|\overline{Z}_s^{i+1}-\overline{Z}_s^{i}|)ds\\
&&+\int_t^Tg(s,\overline{Y}_s^{i+1},\overline{Z}_s^{i+1})dB_s-\int_t^T\overline{Z}_s^{i+1}dW_s.
\end{eqnarray}
For any positive integer $i$, it is obvious that
BDSDE \ref{eq:18} has a unique adapted solution.
By similar procedures, we get the following result:

\begin{theorem}\label{thm:3.9}
Under the assumptions (H1) and (H3)-(H5), and assuming
$\{(\overline{Y}_t^i$, $\overline{Z}_t^i)_{t\in[0,T]}\}_{i=1}^{\infty}$
are the solutions of BDSDEs (\ref{eq:18}), then
\begin{enumerate}
\item[{\rm (i)}]
$\underline{Y}^0_t\leq\overline{Y}^{i+1}_t\leq\overline{Y}^{i}_t\leq\overline{Y}^0_t
$, P-a.s., $\forall t\in [0,T], i=0,1,\cdots$,
\item[{\rm (ii)}]
$\{(\overline{Y}_t^i,\overline{Z}_t^i)_{t\in[0,T]}\}_{i=1}^{\infty}$
converge in $\in S^2([0,T];\mathbb{R})\times M^2(0,T;\mathbb{R}^d)$
to a limit $(\overline{Y}_t,\overline{Z}_t)_{t\in[0,T]}$, which is
the upper bound of solution of  (\ref{eq:1}), in the sense that, for
any other solution $(Y_t,Z_t)_{t\in[0,T]}$ of BDSDE (\ref{eq:1}),
 we have $\overline{Y}_t \geq Y_t$, P-a.s.,
$\forall t\in [0,T].$
\end{enumerate}
\end{theorem}

\begin{remark}\label{rmk:3.10}
It is uncertain whether the upper bound $(\overline{Y}_t,\overline{Z}_t)_{t\in[0,T]}$
of solution of (\ref{eq:1}) is the solution of (\ref{eq:1}).
\end{remark}

\section{Comparison Theorem}\label{sec:4}

The comparison theorem is an important and effective technique in
the theory of BDSDEs. Shi  et al. \cite{SGL} have proved
the comparison theorem for the solutions of BDSDEs with Lipschitz
coefficients. As an application, they showed the existence
of a solution for one dimensional BDSDEs when the coefficient is
continuous with linear growth. In this section, we generalize
the comparison theorem to the case where the coefficient
is left-Lipschitz in $y$ (may be discontinuous) and Lipschitz in $z$.

\begin{theorem}\label{thm:4.1}
(Comparison theorem)\ Let $(Y_t^i,Z_t^i)_{t\in[0,T]}$ $(i=1,2)$ be
the minimal solutions of the following BDSDEs
\begin{eqnarray}\label{eq:19}
Y^1_t=\xi^1+\int_t^Tf_1(s,Y^1_s,Z^1_s)ds+\int_t^Tg(s,Y^1_s,Z^1_s)dB_s-\int_t^TZ^1_sdW_s,
\end{eqnarray}
\begin{eqnarray}\label{eq:20}
Y^2_t=\xi^2+\int_t^Tf_2(s,Y^2_s,Z^2_s)ds+\int_t^Tg(s,Y^2_s,Z^2_s)dB_s-\int_t^TZ^2_sdW_s,
\end{eqnarray}
respectively, where $f_1, f_2, g$ satisfy (H1) and (H3)-(H5),
$\xi^1,\xi^2 \in L^2(\Omega,{\cal{F}}_T,P)$, and $\xi^1\geq\xi^2$, a.s.,
$f_1(t,y,z)$ $\geq f_2(t,y,z)$, a.s., $\forall (t,y,z)\in [0,T]\times
\mathbb{R}\times \mathbb{R}^d$. Then $Y_t^1\geq Y_t^2$, P-a.s.,
$\forall t\in [0,T].$
\end{theorem}
{\bf Proof.}\  Let $(y_t^i,z_t^i)_{t\in[0,T]}$ $(i=0,1,\cdots)$
be the solutions of the following BDSDEs
\begin{eqnarray}\label{eq:21}
y_t^{i}&=&\xi^2+\int_t^T(f_2(s,y_s^{i-1},
z_s^{i-1})-K(y_s^{i}-y_s^{i-1})- K|z_s^{i}-z_s^{i-1}|)ds\\
&&+\int_t^Tg(s,y_s^{i},z_s^{i})dB_s-\int_t^Tz_s^{i}dW_s,\hspace{0.1cm}i\geq 1,
\end{eqnarray}
\begin{eqnarray}\label{eq:22}
y_t^0=\xi^2+\int_t^T(-K|y_s^0|-K|z_s^0|-K) ds
+\int_t^Tg(s,y_s^0,z_s^0)dB_s -\int_t^Tz_s^0dW_s.
\end{eqnarray}
First, we prove $ Y^{1}_t\geq y^0_t$, P-a.s., $\forall t\in [0,T]$.
From (\ref{eq:19}) and (\ref{eq:22}), we have
\begin{eqnarray*}
Y^{1}_t-y^0_t&=&\xi^1-\xi^2+\int_t^T(f_1(s,Y^1_s,Z^1_s)+K|y_s^0|+K|z_s^0|
+K)ds\\
&&+\int_t^T(g(s,Y_s^1,Z_s^1)-g(s,y_s^0,z_s^0))dB_s-\int_t^T(Z_s^1-z_s^0)dW_s\\
&=&\int_t^T(-K|Y^{1}_s-y^0_s|
-K|Z_s^1-z_s^0|+\theta_s^0)ds\\
&&+\int_t^T(g(s,Y_s^1,Z_s^1)-g(s,y_s^0,z_s^0))dB_s-\int_t^T((Z_s^1-z_s^0)dW_s
\end{eqnarray*}
where \begin{eqnarray*}
\theta_s^0:&=&K|Y^{1}_s-y^0_s|
+K|Z_s^1-z_s^0|+K|y_s^0|+K|z_s^0| +K-f_2(s,Y^2_s,Z^2_s)\\
&\geq&K|Y_s^1|+K|Z_s^1| +K+f_1(s,Y^1_s,Z^1_s)\geq 0,
\end{eqnarray*}
we know $ Y^{1}_t\geq y^0_t$, P-a.s., $\forall t\in [0,T]$.

We assume that $ Y^{1}_t\geq y^{i}_t$, P-a.s., $\forall t\in [0,T]$,
from (\ref{eq:19}) and (\ref{eq:21}), we have
\begin{eqnarray*}
&&Y^{1}_t-y^{i+1}_t\\
&=&\xi^1-\xi^2+\int_t^T(f_1(s,Y^1_s,Z^1_s)-f_2(s,y_s^{i},
z_s^{i})+K(y_s^{i+1}-y_s^{i})\\
&&+ K|z_s^{i+1}-z_s^{i+1}|)ds\\
&&+\int_t^T(g(s,Y_s^1,Z_s^1)-g(s,y_s^{i+1},z_s^{i+1}))dB_s-\int_t^T(Z_s^1-z_s^{i+1})dW_s\\
&=&\xi^1-\xi^2+\int_t^T(-K(Y^{1}_s-y^{i+1}_s)
-K|Z_s^1-z_s^{i+1}|+\theta_s^{i+1})ds\\
&&+\int_t^T(g(s,Y_s^1,Z_s^1)-g(s,y_s^{i+1},z_s^{i+1}))dB_s-\int_t^T(Z_s^1-z_s^{i+1})dW_s
\end{eqnarray*}
where
\begin{eqnarray*}
\theta_s^{i+1}&:=&K(Y^{1}_s-y^{i+1}_s)
+K|Z_s^1-z_s^{i+1}|+f_1(s,Y^1_s,Z^1_s)-f_2(s,y_s^{i},
z_s^{i})\\
&&+K(y_s^{i+1}-y_s^{i})+ K|z_s^{i+1}-z_s^{i+1}|\\
&\geq& f_1(s,Y^1_s,Z^1_s)-f_2(s,y_s^{i},
z_s^{i})+K(Y^{1}_s-y^{i}_s)+K|Z_s^1-z_s^{i}|\geq
0,
\end{eqnarray*}
then, for any positive integer $i$,we have $Y^{1}_t\geq y^{i}_t$,
P-a.s., $\forall t\in [0,T]$. From Theorem \ref{thm:3.8}, we have
$\{(y_t^i,z_t^i)_{t\in[0,T]}\}_{i=0}^{\infty}$ converge to
$\{(Y_t^2,Z_t^2)_{t\in[0,T]}\}_{i=0}^{\infty}$. \quad$\Box$

\begin{remark}\label{rmk:4.2}
we replace (H1) by
\begin{enumerate}
\item[(H6)]
$f(t, \cdot, z)$ is right-continuous, and $f(t, y, \cdot)$ is
Lipschitz continuous.
\end{enumerate}
We can deduce
\begin{enumerate}
\item[(i)]
Under the assumptions (H2)-(H4) and (H6),
the same result in Theorem \ref{thm:3.1} holds true.
\item[(ii)]
Under the assumptions (H3)-(H6), we can prove the existence and comparison
results of the maximal solutions of BDSDE (\ref{eq:1}).
\end{enumerate}
\end{remark}

\end{document}